# DISCUSSION OF "EQUI-ENERGY SAMPLER" BY KOU, ZHOU AND WONG

By Ying Nian Wu and Song-Chun Zhu

*University of California, Los Angeles*

We congratulate Kou, Zhou and Wong for making a fundamental contribution to MCMC. Our discussion consists of two parts. First, we ask several questions about the EE-sampler. Then we review a data-driven MCMC scheme for solving computer vision problems.

**1. Questions.** To simplify the language, we use $\pi(x)$ to denote a distribution we want to sample from, and $q(x)$ to denote a distribution at a higher temperature (with energy truncation). The distributions $\pi(x)$ and $q(x)$ can be understood as two consecutive levels in the EE-sampler. Suppose we have obtained a sample from $q(x)$ by running a Markov chain (with a burn-in period), and let us call the sampled states $q$-states. Suppose we have also formed the energy rings by grouping these $q$-states. Now consider sampling from $\pi(x)$ by the EE-sampler.

1. In the jump step, can we make the chain jump to a $q$-state outside the intended energy ring? For example, can we simply propose to jump to any random $q$-states, as if performing the Metropolized independent sampler [2] with $q(x)$ as the proposal distribution? Without restricting the jump to the intended energy ring, it is possible that the chain jumps to a $q$-state of a higher energy level than the current state, but it is also possible that it lands on a lower-energy $q$-state.

If the energy of the current state is very low, we may not have any $q$-states in the corresponding energy ring to make the EE jump. But if we do not restrict the chain to the current energy ring, it may jump to a higher-energy $q$-state and escape the local mode.

The reason we ask this question is that the power of the EE-sampler seems to come from its reuse of the $q$-states, or the long memory of the chain, instead of the EE feature.

2. Can we go up the distribution ladder in a serial fashion? In the EE-sampler, when sampling $\pi(x)$, the chain that samples $q(x)$ keeps running. Can we run a Markov chain toward $q(x)$ long enough and then completely stop it, before going up to sample $\pi(x)$? What is the practical advantage of implementing the sampler in a parallel fashion, or is it just for proving theoretical convergence?

3. About the proof of Theorem 2, can one justify the algorithm by the Metropolized independent sampler [2], where the proposal distribution is $q(x)$ truncated to the current energy range? Of course, there can be a reversibility issue. But in the limit this may not be a problem. In the authors' proof, they also seem to take such a limit. What extra theoretical insights or rigor can be gained from this proof?

4. Can one obtain theoretical results on the rate of convergence? To simplify the situation, consider a Markov chain with a mixture of two moves. One is the regular local MH move. The other is to propose to jump from $x$ to $y$ with $y \sim q$, according to the Metropolized independent sampler. Liu [2] proves that the second largest eigenvalue of the transition kernel of the Metropolized independent sampler is $1 - \min_x \pi(x)/q(x)$. At first sight, this is a discouraging result: even if $q(x)$ captures all the major modes of $\pi(x)$ and takes care of the global structure, the chain can still converge very slowly, because in the surrounding tail areas of the modes, the ratio $\pi(x)/q(x)$ may be very small. In other words, the rate of convergence can be decided by high-energy $x$ that are not important. However, we can regard $q(x)$ as a low-resolution approximation to $\pi(x)$, so we should consider the convergence on a coarsened grid of the state space, where the probability on a coarsened grid point is the sum or integral of probabilities on the underlying finer grid points, so the minimum probability ratio between coarsened $\pi$ and $q$ may not be very small. The lack of resolution in the above scheme is taken care of by the local MH move. So the two types of moves complement each other to take care of things at two different scales. This seems also the case with the more sophisticated EE sampler.

**2. Data-driven MCMC and Swendsen–Wang cut.** Similar to EE-sampler, making large jumps to escape local modes is also the motivation for the data-driven (DD) MCMC scheme of Tu, Chen, Yuille and Zhu [3] for solving computer vision problems.

Let **I** be the observed image data defined on a lattice $\Omega$, and let $W$ be an interpretation of **I** in terms of what is where. One simple example is image segmentation: we want to group pixels into different regions, where the pixel intensities in each region can be described by a coherent generative model. For instance, Figures 1 and 2 show two examples, where the left image is the observed one, and the right image displays the boundaries of the segmented regions. Here $W$ consists of the labels of all the pixels: $W = (W_i, i \in \Omega)$, so that $W_i = l$ if pixel $i$ belongs to region $l \in \{1, \ldots, L\}$, where $L$ is the total number of regions.

In a Bayesian formulation, we have a generative model: $W \sim p(W)$ and $[\mathbf{I}|W] \sim p(\mathbf{I}|W)$. Then image interpretation amounts to sampling from the posterior $p(W|\mathbf{I})$. For the image segmentation problem, the prior $p(W)$ can be something like the Potts model, which encourages identical labels for neighboring pixels. The model $p(\mathbf{I}|W)$ can be such that in each region the pixel values follow a two-dimensional low-order polynomial function plus i.i.d. noise.

To sample $p(W|\mathbf{I})$, one may use a random-scan Gibbs sampler to flip the label of one pixel at a time. However, such local moves can be easily trapped in local modes. A DD-MCMC scheme is to cluster pixels based on local image features, and flip all the pixels in one cluster together.

Specifically, for two neighboring pixels $i$ and $j$, let $p_{i,j} = P(W_i = W_j | F_{i,j}(\mathbf{I}))$, where $F_{i,j}(\mathbf{I})$ is a similarity measure, for example, $F_{i,j}(\mathbf{I}) = |\mathbf{I}_i - \mathbf{I}_j|$. In principle, this conditional probability can be learned from training images with known segmentations. Then for each pair of neighboring pixels $(i, j)$ that belong to the same region under the current state $W = A$, we connect $i$ and $j$ with probability $p_{i,j}$. This gives rise to a number of clusters, where each cluster is a connected graph of pixels. We then randomly pick a cluster $V_0$ (see Figure 3), and assign a single label $l$ to all the pixels in $V_0$ with probability $q_l$. One can design $q_l$ so that the move is always accepted, very much like the Gibbs sampler. This is the basic idea of the Swendsen–Wang cut algorithm of Barbu and Zhu [1], which is a special case of DD-MCMC [4]. The algorithm is very efficient. Figures 1 and 2 show two examples where the results are obtained in seconds, thousands of times faster than the single-site Gibbs sampler.

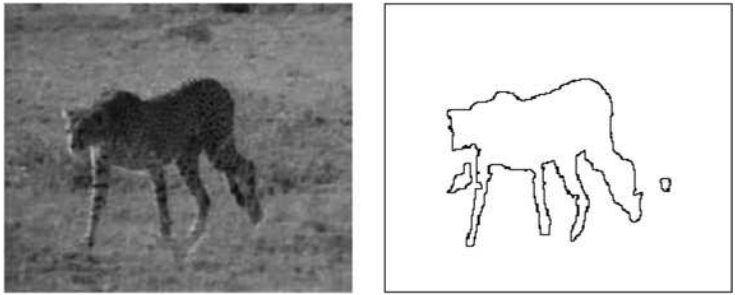

Fig. 1. *Image segmentation.*

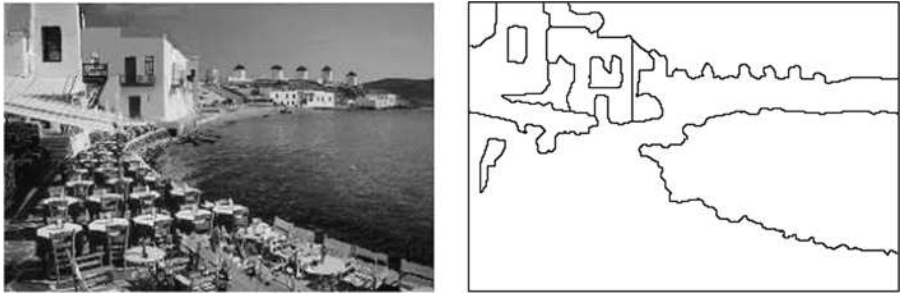

Fig. 2. *Image segmentation.*

Figure 4 illustrates the general situation for DD-MCMC. Part (a) illustrates the model-based inference, where the top–down generative model $p(W)$ and $p(\mathbf{I}|W)$ is explicitly specified. The posterior $p(W|\mathbf{I})$ is implicit and may require MCMC sampling. Part (b) illustrates the bottom-up operations, where some aspects of $W$ can be explicitly calculated based on some simple image features $\{F_k(\mathbf{I})\}$, without an explicit generative model. The bottom-up approach may not give a consistent and accurate full interpretation $W$, but it can be employed to design efficient moves for sampling the posterior $p(W|\mathbf{I})$ in the top–down approach. If vision is a bag of bottom-up tricks, then DD-MCMC provides a principled scheme to bag these tricks. The recent work of Tu, Chen, Yuille and Zhu [3] also incorporates boosting into this MCMC scheme.

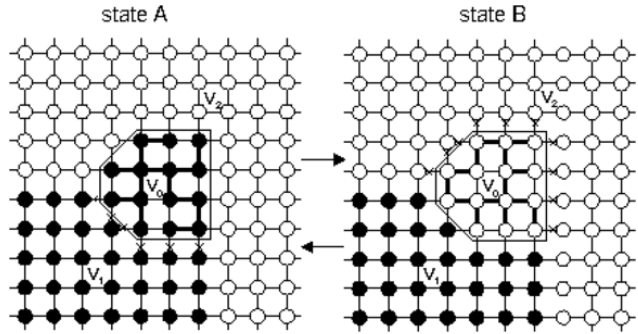

Fig. 3. *Swendsen–Wang cut.*

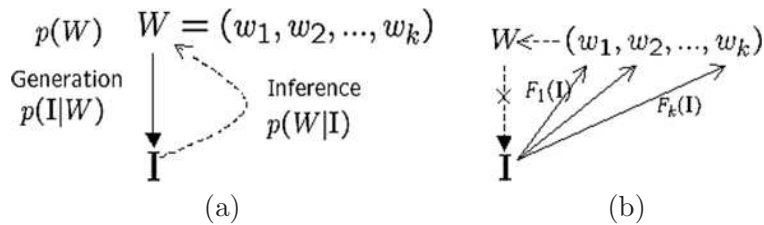

Fig. 4. (a) *Top–down approach;* (b) *bottom-up approach.*

DEPARTMENTS OF STATISTICS
 AND COMPUTER SCIENCE
UNIVERSITY OF CALIFORNIA
LOS ANGELES, CALIFORNIA 90095
USA
E-MAIL: ywu@stat.ucla.edu
          sczhu@stat.ucla.edu